\documentclass[10pt, a4paper, oneside]{article}
\usepackage{amsmath, amssymb, amsthm, amsfonts, amsxtra, latexsym, amscd,
pb-diagram,  graphics}
\usepackage [all]{xy}
\theoremstyle{plain}

\newcommand{\tx}{\otimes }

\newcommand{\ri}{\rightarrow }

\newcommand{\Fm}{\widetilde{F}}

\newcommand{\Gm}{\widetilde{G}}

\theoremstyle{definition}

\begin{document}

\input diagxy

\centerline{\large \bf On monoidal functors between (braided) Gr-categories}
\vspace{0.5cm}
\centerline{\bf Nguyen Tien Quang}
\centerline{\it Hanoi National University of Education, Department of Mathematics,}
\centerline{\it E-mail: nguyenquang272002@gmail.com}
\vspace{0.25cm}
\centerline{\bf Nguyen Thu Thuy}
\centerline{\it Academy of Finance, Science Faculty,}
\centerline{\it Hanoi, Vietnam}
\centerline{\it E-mail: ntthuy11@gmail.com}
\vspace{0.25cm}
\centerline{\bf Pham Thi Cuc}
\centerline{\it Hong Duc University, Natural Science Department,}
\centerline{\it Thanhhoa, Vietnam}
\centerline{\it E-mail: cucphamhd@gmail.com}
\pagestyle{myheadings} 
\markboth{P. T. Cuc, N. T. Quang, N. T. Thuy}{On monoidal
functors between (braided) Gr-categories}


\setcounter{tocdepth}{1}
\vspace{0.5cm}
\noindent{\bf Abstract.} In this paper, we state and prove precise theorems on the classification of the
category of (braided) categorical groups and their (braided)
monoidal functors, and some applications obtained from the basic
studies on monoidal functors between categorical groups.
\\
{\bf 1\ \ Introduction and Preliminaries}\\
\\
Monoidal categories (symmetric monoidal categories) can be ``refined'' to become  {\it
categories with (abelian) group structure} if the notion of  {\it invertible
objects} is added (see [7], [13]). Then, if the underlying category
is a {\it groupoid,} we have the notion of {\it (symmetric) categorical groups}, or {\it Gr-categories} ({\it Picard
categories}) (see [14]). The structure of  Gr-categories and Picard
categories was deeply dealt with by H. X. Sinh in [14].
{\it Braided categorical groups} were originally considered by A. Joyal and R.
Street in [5] as extensions of Picard categories. The category
$\mathcal {BCG}$ of braided categorical groups and braided monoidal
functors was classified by the category $\mathcal{Q}uad$ of
{\it quadratic functions.} These classification theorems were applied and
extended in works of (braided) graded categorical groups (see [1],
[4]), and they led to many interesting results.  It is interesting to revisit even the most basic theory of monoidal categories for further improvements.

In this paper, we state quite adequate studies on Gr-functors
and use them as a common technique to state classification theorems
of the category of categorical groups and the category of braided
categorical groups. Applications motivate these basic studies.

The plan of this paper is, briefly, as follows. In the first section we
recall the construction of a Gr-category of the type $(\Pi,A,h)$, a
reduction of an arbitrary Gr-category.

In the second section, we prove that each Gr-functor between reduced
Gr-categories is one of the type $(\varphi,f)$. Then, we introduce
the notion of the obstruction of a functor of type $(\varphi,f)$,
and cohomological classify these functors.

The next section is devoted to showing the precise theorem of the category of
Gr-categories and Gr-functors, a more complete  version of the
Classification Theorem of H. X. Sinh.

As a consequence of Section 3, it is appropriate to have a different
version of the Classification Theorem of A. Joyal and R. Street for
the category of braided Gr-categories and braided Gr-functors, and
this is the goal of Section 4.

The following section is delicated to the first application  of the obstruction
theory of Gr-functors. We construct the Gr-category of an {\it
abstract kernel} as an example for general theory. This leads to an
interesting consequence: a Gr-category can be transformed into a
strict one (H. X. Sinh proved this result in a completely different way
[15]).

In the last section, we focus on using the Gr-category of an
abstract kernel to classify group extensions by means of
Gr-functors. Then we obtain well-known results on the group
extension problem.

We would stress that Theorem 5 is used in the  method of factor sets to
introduce a new proof of the Classification Theorem for
 graded Gr-categories (see [11]), and a version of Theorem 5 for Ann-functors (see [9]) is used
  to classify Ann-functors thanks to the Mac Lane
cohomology.

 Let us start by some elementary concepts of monoidal
 categories.\\
\indent {\it A monoidal category } $(\mathbb G ,\tx, I, \mathbf{a},
\mathbf{l}, \mathbf{r})$  is a category $\mathbb G$ together with a
tensor product $\otimes: \mathbb G \times \mathbb G\rightarrow
\mathbb G$; and an object $I,$ called the \emph{unit object} of the
category and the natural isomorphisms
\begin{gather*}
\mathbf{a}_{X, Y, Z}:X\tx(Y\tx Z)\rightarrow (X\tx Y)\tx Z,\\
\mathbf{l}_X:I\tx X\ri X\ ,\ \mathbf{r}_X:X\tx I\ri X,
\end{gather*}
which are, respectively, called \emph{associativity, left} and
\emph{right unit constraints}. These constraints  satisfy the
pentagon axiom
\begin{equation*}
(\mathbf{a}_{X, Y, Z}\tx id_T)\  \mathbf{a}_{X,Y\tx Z,T}\ (id_X\tx
\mathbf{a}_{Y,Z,T})=\mathbf{a}_{X\tx Y,Z,T}\ \mathbf{a}_{X,Y,Z\tx
T},
\end{equation*}
and the triangle axiom
$$
id_X\tx \mathbf{l}_Y=(\mathbf{r}_X\tx
id_Y)\mathbf{a}_{X,I,Y}.
$$
\indent A monoidal category is \emph{strict} if the associativity
constraint $\mathbf{a}$
and the unit constraints  $\mathbf{l},\mathbf{r}$ are all identities.

\indent Let $\mathbb G=(\mathbb G, \tx, I, \mathbf{a}, \mathbf{l},
\mathbf{r})$ and $\mathbb G'=(\mathbb G', \tx, I', \mathbf{a}',
\mathbf{l}', \mathbf{r}')$ be monoidal categories. \emph{A monoidal
functor} from $\mathbb G$ to $\mathbb G'$ is a triplet $(F,
\Fm,F_\ast)$ where  $F:\mathbb G \ri\mathbb G'$  is a functor,
$F_\ast:I'\ri FI$ is an isomorphism, and natural isomorphisms
$$\Fm_{X,Y}:FX\tx FY\ri F(X\tx Y)$$
 such that
\begin{gather*}
F(\mathbf{a}_{X,Y,Z})\circ \widetilde{F}_{X,YZ}\circ (FX\otimes \widetilde{F}_{Y,Z})=\widetilde{F}_{X\otimes Y,Z}\circ (\widetilde{F}_{X,Y}\otimes FZ)\circ \mathbf{a'}_{FX,FY,FZ},\\
\mathbf{r}'_{FX}=F(\mathbf{r}_X)\circ \widetilde{F}_{X,I}\circ
(id\otimes F_\ast):FX\otimes I'\rightarrow FX,\\
\mathbf{l}'_{FX}=F(\mathbf{l}_X)\circ \widetilde{F}_{I,X}\circ (
F_\ast\otimes id):I'\otimes FX \rightarrow FX.
\end{gather*}

\indent \emph{A natural monoidal equivalence} or a {\it homotopy}
$\alpha:(F,\Fm,F_\ast)\ri  (G,\Gm, G_\ast)$ between monoidal
functors from $\mathbb G$ to $\mathbb G'$ is a natural isomorphism
$\alpha: F\rightarrow G,$ such that
\begin{equation*}
G_\ast=\alpha_I\circ F_\ast, 
\end{equation*}
and
\begin{equation*}
\alpha_{X\otimes Y}\circ\widetilde{F}_{X,Y}=\widetilde{G}_{X,Y}\circ (\alpha_X\otimes\alpha_Y). 
\end{equation*}
\indent 
\emph{A monoidal equivalence} between monoidal categories is a
monoidal functor $F:\mathbb G \rightarrow \mathbb G',$ such that
there exists a monoidal functor $G:\mathbb G'\rightarrow \mathbb G$
and homotopies $\alpha: G.F\rightarrow id_{\mathbb G}$ and $
\beta:F.G\rightarrow id_{\mathbb G'}$. $(F, \Fm, F_\ast)$ is a
monoidal equivalence if and only if $F$ is an equivalence.


A {\it Gr-category} is a monoidal category, where every object is
invertible and every morphism is isomorphic. If $(F, \Fm, F_\ast)$
is a monoidal functor between Gr-categories, it is called a
{\it Gr-functor.} Then the isomorphism $ F_\ast:I'\ri FI$ can be
deduced from $F$ and $\Fm$.

 Let us recall some well-known results
(see [14]). Each Gr-category $\mathbb G$ is equivalent to a
Gr-category of the type $(\Pi, A),$ which can be described as
follows. The set $\pi_0\mathbb G$ of iso-classes of objects of
$\mathbb G$ is a group with the operation induced by the tensor
product in $\mathbb G,$ and the set $\pi_1\mathbb G$ of
automorphisms of the unit object $I$ is an abelian group with the
operation, denoted by +, induced by the composition of morphisms.
Moreover, $\pi_1\mathbb G$ is a $\pi_0\mathbb G$-module with the
action
$$su=\gamma_{X}^{-1}\delta_{X}(u),\ X\in s,\ s\in \pi_0\mathbb G,\ u\in \pi_1\mathbb G,$$
where $\delta_X,\ \gamma_X$ are defined by the following commutative
diagrams
$$\begin{CD}
X @>\gamma_X(u) >> X \\
@A {\bf l}_X AA     @AA {\bf l}_X A \\
I \tx X @>u \tx id >> I \tx X
\end{CD}
\qquad \qquad\qquad
\begin{CD}
X @>\delta_X(u) >> X \\
@A {\bf r}_X AA     @AA {\bf r}_X A \\
X \tx I @>id \tx u >> X \tx I
\end{CD}$$

\indent The reduced Gr-category $S_{\mathbb G}$ of a Gr-category
$\mathbb G$ is a category whose objects are the elements of
$\pi_0\mathbb G$ and whose morphisms are automorphisms $(s,u):s\ri
s,$ where $s \in \pi_0\mathbb G,\ u \in \pi_1\mathbb G$. The
composition of two morphisms is induced by the addition in
$\pi_1\mathbb G$
$$(s,u).(s,v)=(s,u+v).$$
The category $S_{\mathbb G}$  is equivalent to $\mathbb G$ by
canonical equivalences constructed as follows. For each
$s=[X] \in \pi_0\mathbb G,$ choose a representative $X_s\in \mathbb
G;$ and for each $X \in s,$ choose an isomorphism $i_X:\
X_s\rightarrow X$ such that $i_{X_s}=id_{X_s}$. For the set of representatives, we obtain two functors
\begin{equation}\label{ct1}
\begin{cases}
G:\mathbb G\ri S_{\mathbb G} ,\\
G(X)=[ X]=s,\\
G(X \stackrel {f}{\ri}Y)=(s,\gamma_{X_s}^{-1}(i_{Y}^{-1}fi_X)),
\end{cases}\qquad\qquad
\begin{cases}
H: S_{\mathbb G} \ri  \mathbb G,\\
H(s)=X_s,\\
H(s,u)=\gamma_{X_s}(u).
\end{cases}
\end{equation}
\indent Two functors $G$ and $H$ are categorical equivalences by
natural transformations
$$\alpha=(i_{X}): HG\cong id_{\mathbb G};\qquad\ \beta=id:GH\cong id_{S_{\mathbb G}}.$$
They are called \emph{canonical equivalences}.

With the structure transport (see [13], [14]) by the quadruple $(G,
H, \alpha, \beta),$ $S_{\mathbb G}$  becomes a Gr-category together
with the following operation:
\begin{gather*}
s\tx t=s.t,\quad s,t\in\pi_0 \mathbb G,\\
(s,u)\tx(t,v)=(st,u+sv),\quad u,v\in\pi_1 \mathbb G.
\end{gather*}
\noindent The set of representatives $(X_s, i_X)$ is called a
\emph{stick} of the Gr-category $\mathbb G$ for
$$X_1=I,\;i_{I\otimes X_s}=\mathbf{l}_{X_s},\;i_{X_s\otimes
I}=\mathbf{r}_{X_s}.$$ The unit constraints of the Gr-category
$S_\mathbb G$ are therefore strict, and its associativity constraint
  is a normalized 3-cocycle $h\in Z^{3}(\pi_0\mathbb G,\pi_1\mathbb G)$.
Moreover, the equivalences $G,\ H$ become Gr-equivalences together
with natural isomorphisms
\begin{equation}\label{ct2}
\widetilde{G}_{A,B}=G(i_A\tx i_B),\
 \widetilde{H}_{s,t}=i_{X_s\tx X_t}^{-1}.
 \end{equation}
The Gr-category $S_{\mathbb G}$  is called a \emph{reduction} of the
Gr-category $\mathbb G.$ $S_{\mathbb G}$  is said to be of the
{\it type} $(\Pi, A, h)$ or simply the {\it type} $(\Pi,A)$ if
$\pi_0\mathbb G, \pi_1\mathbb G$ are, respectively, replaced with
the group $\Pi$ and the $\Pi$-module $A.$\\
\\
{\bf 2\ \ Classification of Gr-functors of the type $(\varphi,f)$}\\
\\
In this section, we show that each  Gr-functor   $(F,
\Fm):\mathbb G\ri\mathbb G^{'}$ induces a  Gr-functor $S_F$ between
their reduced Gr-categories. This allows us to study the problem of
the existence of  Gr-functors and classify them on  Gr-categories of
the type  $(\Pi, A).$ The following proposition is mentioned in many
works related to categorical groups.

\noindent {\bf Proposition 2.1.} [14] {\it Let $(F, \Fm): \mathbb G\ri\mathbb G{'}$ be a Gr-functor. Then,
$(F,\Fm)$ induces a pair of group homomorphisms
\begin{gather*}
 F_0: \pi_{0}\mathbb G\rightarrow \pi_{0}\mathbb G^{'}, \ \ [X]\mapsto [FX],\\
F_1:\pi_{1}\mathbb G\rightarrow \pi_{1}\mathbb G^{'}, \ \ u\mapsto \gamma^{-1}_{FI}(Fu)
\end{gather*}
satisfying $F_1(su)= F_0(s)F_1(u)$.}

Our first result is to strengthen Proposition 2.1 by Proposition 2.4 as
asserting that each Gr-functor $(F, \Fm): \mathbb G\ri\mathbb G{'}$
induces a Gr-functor $S_{\mathbb G}\ri S_{\mathbb G}{'}.$ We need
two following  lemmas

\noindent {\bf Lemma 2.2.} {\it Let $\mathbb G, \mathbb G'$ be two $\otimes$-categories with,
 respectively, unit constraints $(I,{\bf l},{\bf r})$, $(I',{\bf l'},{\bf r'})$, and $(F,\Fm,F_{\ast}):\mathbb G
\ri \mathbb G'$  be an $\otimes$-functor which is compatible with
  the unit constraints. Then, the following diagram commutes:
\[\begin{diagram}
\node{FI}\arrow{e,t}{\gamma_{_{FI}}(u)}
\node{FI}\\
\node{I'}\arrow{n,l}{F_{\ast}}\arrow{e,t}{u}
\node{I'}\arrow{n,r}{F_{\ast}}
\end{diagram}\]
It follows that
$$\gamma^{-1}_{FI}(Fu)=F^{-1}_{\ast}F(u)F_{\ast},$$
\noindent i.e., the notions of the map $F_1$ in [1] and in
Proposition 1 coincide.}

\begin{proof}  Clearly, $\gamma_{I'}(u)=u$. Moreover, the family
$(\gamma_{X'}(u)), \ X' \in \mathbb{G'},$ is an endomorphism of the identity
functor $id_{\mathbb{G'}}.$ So the above diagram commutes.

The final conclusion is deduced from the above commutative diagram,
when $u$ is replaced by $\gamma^{-1}_{FI}(Fu).$
\end{proof}
\noindent {\bf Lemma 2.3.} {\it If the hypothesis of Lemma 2.2 holds, we have}
 $$F\gamma_{_X}(u) = \gamma_{_{FX}}(\gamma^{-1}_{_{FI}}Fu).$$
\begin{proof}
Consider the following diagram
$$\bfig
\square |almb|<700,500>[I'\tx FX`FI \tx FX` I'\tx FX ` FI \tx
FX;F_{\ast} \tx id `\gamma^{-1}_{FI}Fu \tx id `Fu \tx id`F_{\ast}
\tx id] \square (700,0)|bmma|/>``>`>/<700,500>[FI\tx FX`F(I \tx X)`
FI \tx FX`F(I \tx X);\Fm ``F(u\tx id) `\Fm] \square
(1400,0)|amrb|/>``>`>/<700,500>[F(I \tx X)`FX`F(I \tx
X)`FX;F(\mathbf{l}_X) ``F\gamma_X(u)`F(\mathbf{l}_X)] \morphism
(0,500)/-/ <0,250>[I'\tx FX`;] \morphism (0,750)/-/
<2100,0>[`;\mathbf{l'}] \morphism (2100,750)/->/ <0,-250>[`FX;]
\morphism (0,0)/-/ <0,-250>[I'\tx FX`;] \morphism (0,-250)|b|/-/
<2100,0>[`;\mathbf{l'}] \morphism (2100,-250)/->/ <0,250>[`FX;]
\place(350,250)[(1)]\place(1050,250)[(2)]\place(1750,250)[(3)]\place(1050,650)[(5)]
\place(1050,-150)[(4)] \efig $$ In this diagram, the regions (4) and
(5) commute thanks to the compatibility of the functor
$(F,\widetilde{F})$ with the unit constraints. The region (3)
commutes due to the definition of $\gamma_X$ (with image through
$F$), the region (1) commutes by Lemma 2.2. The commutativity of the region (2) follows from the naturality of the isomorphism $\widetilde{F}.$
Therefore, the outer perimeter commutes, i.e.,
$F\gamma_X(u)=\gamma_{FX}\big(\gamma^{-1}_{FI}Fu\big).$
\end{proof}
{\it Remark on notations:} Hereafter, if there is no explanation, $\mathbb S,\mathbb S'$ refer to Gr-categories $(\Pi, A,h),(\Pi',A',h').$

A functor $F: \mathbb S\rightarrow \mathbb S'$ is
called a functor of the {\it type} $(\varphi,f)$ if
\begin{equation*}
F(x)=\varphi(x),\ \ F(x,u)=(\varphi(x), f(u)),
\end{equation*}
where $\varphi:\Pi\rightarrow \Pi'$, $f:A\rightarrow A'$ is a pair
of group homomorphisms satisfying $f(xa)=\varphi (x)f(a)$ for $x\in
\Pi, a\in A.$

\noindent {\bf Proposition 2.4.} {\it Each Gr-functor $(F, \Fm): \mathbb G\rightarrow \mathbb G'$ induces a
Gr-functor $S_F:S_{\mathbb G}\rightarrow S_{\mathbb G'}$ of the type
$(\varphi, f),$
 with $\varphi=F_0, f=F_1$. Moreover, $S_F=G' F H,$
where $H, G'$ are canonical equivalences.}


\begin{proof} Let $K = G'F H$ be the composition functor.  One can verify that $K(s)= F_0(s)$, for $s\in
\pi_0\mathbb G$. We now prove that $K(s,u)= (F_0s, F_1u)$ for each
morphism $u: I \ri I.$ We have
$$K(s,u) = G'FH(s,u) = G'(F \gamma_{X_s}(u)).$$
 \indent Since $H'G' \simeq id_{\mathbb{G'}},$ by the natural
 equivalence $\beta = (i'_{X'}),$ the following diagram commutes (note that $X_s' = H'G' FX_s$):
$$\begin{CD}
X'_{s'}  @>i'>>  FX_s \\
@V H'G'F \gamma_{X_s}(u)VV   @VV F \gamma_{X_s}(u) V \\
X'_{s'} @>i'>> FX_s
\end{CD}$$
\indent According to Lemma 2.3, we have
$$ F\gamma_{X_s}(u)=\gamma_{FX_s}(\gamma^{-1}_{FI}Fu).$$
\indent Besides, since the family $(\gamma_{X'})$ is a natural
equivalence of the identity functor $id_{\mathbb{G}'}$, the
following diagram commutes:
\[\begin{CD}
X'_{s'}  @>i'>>  FX_s \\
@V \gamma_{X_s'}(\gamma^{-1}_{FI}Fu)VV   @VV \gamma_{FX_s}(\gamma^{-1}_{FI}Fu) V \\
X'_{s'}  @>i'>>  FX_s
\end{CD}\]
Hence,  $H'G'F \gamma_{X_s}(u)=
\gamma_{X'_s}\big(\gamma^{-1}_{FI}Fu\big)$. By the definition of
$H',$ we have
$$G'F\gamma_{X_s}(u)=(F_0s,\gamma^{-1}_{FI}Fu)=(F_0s,F_1(u)).$$
 This means $K=S_F$.
\end{proof}

We now describe Gr-functors on Gr-categories of the type $(\Pi, A)$.

\noindent {\bf Theorem 2.5.} {\it Any  $(F,\Fm):\mathbb S\ri \mathbb S'$ is a  Gr-functor of the type $ ( \varphi,
f).$}
\begin{proof}
\indent For $x,y\in \Pi$, $\Fm_{x,y}:Fx\tx Fy\ri F(x\tx y)$ is a
morphism in $\mathbb S'$. It follows that $Fx.Fy=F(xy)$. So if one
sets $\varphi(x)=Fx$,  $\varphi:\Pi\rightarrow\Pi'$ is a group
homomorphism.

\indent  We write $F(x,a)=(\varphi(x),f_x(a))$. Since $F$ is a
functor, we have
$$F((x,a).(x,b))=F(x,a).F(x,b).$$
It follows that
\begin{equation*}
 f_x(a+b)=f_x(a)+f_x(b).
\end{equation*}
Therefore, $f_x:A\ri A'$ is a group homomorphism for  each $x\in
\Pi$. Besides, since $(F,\widetilde{F})$ is an $\otimes$-functor,
the following diagram commutes
$$\begin{CD}
Fx.Fy @>\Fm >> F(xy)\\
@V Fu\tx Fv VV   @ VV F(u\tx v) V \\
Fx.Fy @ >\Fm >> F(xy)
\end{CD}$$
for all   $u=(x,a),\ v=(y,b)$. Hence, we have
\begin{gather*}
\ \ F(u\tx v)=Fu\tx Fv\\
\qquad \Leftrightarrow\   f_{xy}(a+xb)=f_x(a)+\varphi(x).f_y(b)
\end{gather*}
\begin{equation}\label{ct3}
\Leftrightarrow\  f_{xy}(a)+f_{xy}(xb)=f_x(a)+\varphi(x).f_y(b).
\end{equation}
In the relation (\ref{ct3}), let $x=1,$ we obtain $f_y(a)=f_1(a)$. Thus,
$f_y=f_1$ for all $y\in \Pi.$ Write $f_y=f$ and use (\ref{ct3}), we obtain
$f(xb)=\varphi(x).f(b).$
\end{proof}

Note that if $\Pi'$-module $A'$ is regarded as a $\Pi$-module by the
action  $xa'=\varphi(x).a'$, then $f:A\ri A'$ is a homomorphism of
$\Pi$-modules. Since $$\Fm_{x,y}=(F(xy), g_F(x,y)):Fx.Fy\ri F(xy),$$
where $g_F:\Pi^2\ri A'$  is a function, it is said that $g_F$ is {\it
associated}  with $\widetilde{F}$. The compatibility of $(F, \Fm)$
with the associativity constraint leads to the relation:
$$\varphi^\ast h'-f_\ast h=\partial(g_F),$$
where
\begin{gather*}
 (f_{\ast} h)(x,y,z)=f( h(x,y,z)),\\
(\varphi^{\ast} h')(x,y,z)= h'(\varphi x, \varphi y,\varphi z).
\end{gather*}
One can see that  two Gr-functors $(F,\widetilde{F}), (F',\widetilde{F}'):\mathbb S\rightarrow \mathbb S'$ are homotopic if and
only if $F'=F$, i.e., they are of the same type $(\varphi, f)$, and
there is a function $t:\Pi\ri A'$ such that $g_F'=g_F+\partial t.$

 We refer to
$$\text{Hom}_{(\varphi, f)}[\mathbb S, \mathbb S'].$$
as the set of homotopy classes of Gr-functors of the type $(\varphi,
f).$

 In order to find the sufficient condition to make a functor of the
type $(\varphi,f)$ become a Gr-functor, we introduce the notion of
{\it the obstruction} like in the case of Ann-functors (see [10]). If $h,h'$ are,respectively, associativity constraints of Gr-categories $\mathbb S, \mathbb S'$ and $F:\mathbb S\rightarrow \mathbb S'$ is a functor of the type
$(\varphi,f),$ then the function
 \begin{equation}\label{ct4}
 k=\varphi^{\ast} h'-f_{\ast} h
 \end{equation} is called \emph{an obstruction}
of $F.$ 

Keeping in mind that $\mathbb S=(\Pi, A,h),\mathbb S'=(\Pi',A',h'),$ we move on the following theorem

\noindent {\bf Theorem 2.6.} {\it The functor  $F:\mathbb S\ri \mathbb S' $  of the type  $(\varphi, f)$ induces  a
   Gr-functor if
and only if its obstruction $[k]=0$
 in  $H^3(\Pi, A')$. Then, there exist bijections:}
\begin{eqnarray}
\text{i) Hom}_{(\varphi, f)}[\mathbb S, \mathbb S']\rightarrow H^2(\Pi, A'),\hspace{6.7cm}
 \end{eqnarray}
ii) $\text{Aut}(F)\rightarrow  Z^1(\Pi, A').$
\begin {proof}  If $(F,\Fm): \mathbb S\ri \mathbb S'$ is a Gr-functor, then  $(F,\Fm)=(\varphi,f,g_F)$,
where
$$\varphi^\ast{ h'}-f_\ast h= \partial (g_F)\in B^{3}(\Pi,A').$$
 Therefore,
$[\varphi^\ast{ h'}]-[f_\ast h]=0$  in $H^3(\Pi, A')$.

\indent Conversely, since $[\varphi^\ast{ h'}]-[f_\ast h]=0$ there
exists a 2-cochain $g\in Z^2(\Pi,A')$ such that ${\varphi^\ast
h'}-f_\ast h= \partial g$. Take $\Fm$ be associated with $g,$  one
can see that  $(F,\Fm)$ is a Gr-functor.


\noindent i) If  $(F,\Fm):\mathbb S\ri \mathbb S'$ is a  Gr-functor, then  $F=(\varphi, f, g_F)$.
   Let $g_F$ be fixed. Now if
$$(K,\widetilde{K}):\mathbb S\ri \mathbb S'$$
 is a Gr-functor of the type   $(\varphi, f)$, then
 $\partial(g_F)={\varphi^* h'}-f_* h= \partial (g_K).$
 It follows that $g_F-g_K$ is a
2-cocycle. Consider the correspondence
\[\Phi: [(K,\widetilde{K})]\mapsto [g_F-g_K]\]
between the set of congruence classes of Gr-functors of the type
$(\varphi, f)$ from $\mathbb S$ to $\mathbb S'$  and the
group $H^2(\Pi, A').$

First, we show that the above correspondence is a map. Indeed,
let
$$(K',\widetilde{K}'):\mathbb S\ri \mathbb S'$$
\noindent be a Gr-functor and $u: K\ri K'$ be a homotopy. Then $K,
K'$ are of the same type $(\varphi, f)$ and $g_{K'}=g_K+\partial t$
where $g_K, g_{K'}$ are, respectively, associated with
$\widetilde{K}, \widetilde{K}',$
 i.e., $[g_F-g_{K'}]=[g_F-g_K]\in H^2(\Pi, A').$

Furthermore, $\Phi$ is an injection.




Finally, we show that the correspondence $\Phi$ is a
surjection. Indeed, assume that $g$ is an arbitrary 2-cocycle. We
have
\[\partial (g_F-g)=\partial g_F-\partial g=\partial g = \varphi^* h'-f_* h.\]
Then, there exists a Gr-functor
$$(K,\widetilde{K}):\mathbb S\ri \mathbb S'$$
of the type $(\varphi,f),$ with a functor isomorphism
$\widetilde{K}=(\bullet,g_F-g).$ So $\Phi$ is a surjection.

\noindent ii)  Let $F=(F,\Fm):\mathbb S\ri \mathbb S'$ be a Gr-functor and  $t \in$ Aut$(F)$.
Then, the equality  $g_F=g_F+\partial t$ implies that  $\partial
t=0$, i.e., $t\in Z^1(\Pi, A')$.
\end{proof}
\noindent {\bf 3\ \ The general case}\\
\\
Let $\mathcal{CG}$ be a category whose objects are  Gr-categories,
and whose morphisms
are monoidal functors between them. 
We determine the category $\bf H^{3}_{Gr}$ whose objects are
triplets $(\Pi,A,[h])$ where  $[h]\in H^3(\Pi,A)$. A morphism $(\varphi,f):(\Pi,A,[h])\ri (\Pi',A',[h'])$ in $\bf H^{3}_{Gr}$  is
a pair $(\varphi,f)$ such that there is a function $g:\Pi^2\ri A'$ so that
$(\varphi,f,g):(\Pi,A,h)\ri (\Pi',A',h')$ is a Gr-functor, i.e.,
$[\varphi^\ast h']=[f_\ast h]\in H^3(\Pi,A')$. The composition in
$\bf H^{3}_{Gr}$ is given by
$$(\varphi',f')\circ (\varphi,f)=(\varphi'\varphi, f'f).$$
One can observe that {\it two Gr-functors $F, F':\mathbb G \ri \mathbb
G'$ are homotopic if and only if $F_i=F_i', i=0,1$ and
$[g_F]=[g_{F'}]$}. Denote the set of homotopy classes of Gr-functors
$\mathbb G\ri\mathbb G'$ which induce the same the pair $(\varphi,
f)$ by
$$\text{Hom}_{(\varphi, f)}[\mathbb G,\mathbb G'].$$
We now state the main result of this section

\noindent {\bf Theorem 3.1.} (The Classification Theorem) {\it There is a classifying functor
\[\begin{matrix}
 d:&\mathcal{CG}&\ri&\bf H^{3}_{Gr},\\
&\mathbb{G}&\mapsto&(\pi_0\mathbb{G},\pi_1\mathbb{G},[h_\mathbb{G}]),\\
&(F,\widetilde{F})&\mapsto&(F_0,F_1)
\end{matrix}\]
which has the following properties:}\\
 i) $dF$ {\it is an isomorphism if and
only if $F$ is an
 equivalence.}\\
ii) $d$ {\it is a surjection on objects.}\\
 iii) $d$ {\it is full, but not faithful. For $(\varphi,f):d\mathbb{G}\rightarrow
d\mathbb{G}'$, there is a bijection}
\begin{equation}\label{ct6}
\overline{d}:\text{Hom}_{(\varphi,f)}[\mathbb{G},\mathbb{G}']\rightarrow
H^2(\pi_0\mathbb{G},\pi_1\mathbb{G}').
\end{equation}

\begin{proof}
In the Gr-category $\mathbb{G}$, for each stick $(X_s,i_X)$ we can
construct a reduced   Gr-category
$(\pi_0\mathbb{G},\pi_1\mathbb{G},h)$. If the choice of the stick is
modified, then the 3-cocycle  $h$ changes to a cohomologous 3-cocycle $h'.$ Therefore,  $\mathbb{G}$ determines
a unique element $[h]\in H^3(\pi_0\mathbb{G},\pi_1\mathbb{G})$. This
shows that $d$ is a map on objects.

For Gr-functors
$$\mathbb{G}\stackrel{F}{\ri}\mathbb{G}'\stackrel{F'}{\ri}\mathbb{G}'',$$
\noindent one can see that $(F'F)_0=F'_0F_0$. Since $(F'F)_\ast$
is the composition
$$I''\stackrel{F'_\ast}{\ri}F'I'\stackrel{F'(F_\ast)}{\ri}F'FI,$$
then for $u\in$ Aut$(I)$ we have
\begin{eqnarray*}
(F'F)_1u&=&(F'F)_\ast^{-1}(F'F)u(F'F)_\ast\\
&=&F_*'^{-1}F'(F_\ast^{-1})F'FuF'(F_\ast)F'_\ast\\
&=&F_\ast'^{-1}F'(F_1u)F'_\ast=F'_1(F_1u).
\end{eqnarray*}
That is
$$d(F'\circ F)=(dF')\circ (dF).$$
Clearly, $d({id_\mathbb{G}})=id_{d\mathbb{G}}$. Therefore, $d$ is a
functor.

\noindent i) According to Proposition 1.1.

\noindent ii) If $(\Pi,A,[h])$ is an object of $\bf{H}^3_{Gr}$, $\mathbb
S=(\Pi,A,{h})$ is a Gr-category of the type $(\Pi,A)$ and obviously
$d\mathbb S=(\Pi,A,[h])$ .

\noindent iii) Let $(\varphi,f)$ be a morphism in
$\text{Hom}_{\bf{H}^3_{Gr}}(d\mathbb{G},d\mathbb{G'})$. Then, there
exists
a function $g:(\pi_0\mathbb{G})^2\ri \pi_1\mathbb{G'}$ such that
$$\varphi^\ast{h}_\mathbb{G'}=f_\ast{h}_\mathbb{G}+\partial g.$$
Hence, by Theorem 2.6,
$$K=(\varphi,f,g):(\pi_0\mathbb{G}, \pi_1\mathbb{G},h_\mathbb{G})\ri
(\pi_0\mathbb{G'}, \pi_1\mathbb{G'},h_\mathbb{G'})$$ is a
Gr-functor. Then, the composition  Gr-functor $F=H' K
G:\mathbb{G}\ri\mathbb{G'}$ induces  $dF=(\varphi,f)$.  This shows
that the functor $d$ is full.

To prove that (\ref{ct6}) is a bijection, we prove the correspondence
\begin{equation}\label{ct7}
\overline{s}:\text{Hom}_{(\varphi,f)}[\mathbb{G},\mathbb{G'}]{\ri}
\text{Hom}_{(\varphi,f)}[S_\mathbb{G},S_{\mathbb{G'}}],
\end{equation}
\centerline{$\ \ [F]\mapsto[S_F]$}
is a bijection.

Clearly, if $F, F':\mathbb G\ri \mathbb G'$ are homotopic, then the
induced  Gr-functors $S_F, S_{F'}:S_{\mathbb G}\rightarrow
S_{\mathbb G'}$ are homotopic. Conversely, if $F, F'$ are
Gr-functors such that $S_F, S_{F'}$
 are homotopic, then the
compositions $E=H'(S_F)G$ and $E'=H'(S_{F'})G$ are homotopic, where
$H', G$ are canonical Gr-equivalences. The Gr-functors $E, E'$ are
respectively homotopic to $F, F'$. Hence, $F$ and $F'$ are
homotopic. This shows that $\overline{s}$ is an injection.



Now, if $K=(\varphi,f,g):S_\mathbb{G}\ri S_{\mathbb{G}'}$ is a
Gr-functor, then the composition
$$F=H' K
G:\mathbb{G}\ri\mathbb{G}'$$ is a Gr-functor satisfying $S_F=K$,
i.e.,  $\overline{s}$  is a surjection. Finally, the bijection (\ref{ct6})
is the composition of (5) and (\ref{ct7}).
\end{proof}


By Theorem 3.1, one can simplify the problem of equivalence
classification of Gr-categories by the one of  Gr-categories with
the same (up to an isomorphism) two first invariants.

 Let $\Pi$ be a group and $A$ be a $\Pi$-module. It is said that the
Gr-category $\mathbb G$ has a {\it pre-stick of the type} $(\Pi, A)$
if there exists a pair of group isomorphisms
\[p: \Pi\ri \pi_0\mathbb G,\qquad q: A\ri \pi_1\mathbb G\]
which is compatible with the action of modules
\[q(su)=p(s)q(u),\]
 where $s\in \Pi, u\in A$. The pair $\epsilon=(p, q)$ is called a
{\it pre-stick of the type} $(\Pi, A)$ to the Gr-category $\mathbb
G$.

  A {\it morphism} between the two Gr-categories $\mathbb G, \mathbb
G'$ whose pre-sticks are of the type $(\Pi, A)$  (with,
respectively, the pre-sticks    $\epsilon=(p, q), \epsilon'=(p',
q')$) is a Gr-functor $(F, \Fm): \mathbb G\ri \mathbb G'$ such that
the following diagrams commute
\[\begin{diagram}
\node{\pi_0\mathbb G}\arrow[2]{e,t}{ F_0} \node[2]{\pi_0\mathbb G'}
\node[3]{\pi_1\mathbb G}\arrow[2]{e,t}{F_1}
\node[2]{\pi_1\mathbb G'}\\
\node[2]{\Pi}\arrow{nw,b}{p}\arrow{ne,b}{p'}
\node[5]{A}\arrow{nw,b}{q}\arrow{ne,b}{q'}
\end{diagram}\]
where $F_0, F_1$ are two homomorphisms induced from $(F,\Fm)$.

Clearly, it follows from the definition  that $F_0, F_1$ are
isomorphisms and therefore $F$ is an equivalence. Let
$$\mathcal{CG}[\Pi, A]$$
denote the set of equivalence classes of Gr-categories whose
pre-sticks are of the type  $(\Pi, A)$. We can prove the
Classification Theorem of H. X. Sinh [15] based on these results
as follows.

\noindent {\bf Theorem 3.2.} (H. X. Sinh) {\it There exists a bijection}
\begin{eqnarray*}
	\Gamma:\mathcal{CG}[\Pi, A]&\ri & H^{3}(\Pi, A),\\
{[\mathbb G]}&\mapsto & q^{-1}_{\ast} p^{\ast} [h_{\mathbb G}].
\end{eqnarray*}
\begin {proof}
By Theorem 3.1, each Gr-category $\mathbb G$ determines
 uniquely an element
$[h_{\mathbb G}]\in H^3(\pi_0\mathbb G, \pi_1\mathbb G)$, and then
determines an element
$$\epsilon [h_{\mathbb G}]=q^{-1}_\ast p^\ast [h_{\mathbb G}]\in H^3(\Pi, A).$$
Now, if $F:\mathbb G \ri \mathbb G'$ is a morphism between two
Gr-categories whose pre-sticks of the type $(\Pi, A)$, then the
induced  Gr-functor $S_F=(\varphi, f, g_F)$  satisfies 
$$\varphi^\ast [h_{\mathbb G'}]=f_\ast [h_\mathbb G].$$
It follows that
$$\epsilon' [h_{\mathbb G'}]=\epsilon [h_{\mathbb G}].$$
This means $\Gamma$ is a map. Moreover, it is an injection. Indeed,
suppose that $\Gamma[\mathbb G]=\Gamma[\mathbb G'],$  we have
$$\epsilon' (h_{\mathbb G'})-\epsilon (h_{\mathbb G})=\partial g.$$
Therefore, there exists a Gr-functor $J$ of the type $(id,id)$ from
$\mathbb{J}=(\Pi,A, \epsilon (h_{\mathbb G}))$ to
$\mathbb{J}'=(\Pi,A,\epsilon' (h_{\mathbb G'}))$. The composition
$$\mathbb G\stackrel{G}{\rightarrow} S_{\mathbb G}\stackrel{\epsilon^{-1}}{\rightarrow}
\mathbb{J}\stackrel{J}{\rightarrow}
\mathbb{J}'\stackrel{\epsilon'}{\rightarrow}S_{\mathbb
G'}\stackrel{H'}{\rightarrow}\mathbb G'$$ implies $[\mathbb
G]=[\mathbb G']$, and $\Gamma$ is an injection.
 Obviously, $\Gamma$ is a surjection.

\end {proof}

\noindent {\bf 4\ \ The case of braided Gr-categories}\\
\\
A Gr-category $\mathbb B$ is called a {\it braided Gr-category} if
there is a {\it braiding} ${\bf c}$, i.e., a natural isomorphism  ${\bf
c}={\bf c}_{X,Y}:X\otimes Y\rightarrow Y \otimes X$, which is
compatible with $\mathbf{a}, \mathbf{l}, \mathbf{r}$ in the sense of
satisfying the following coherence conditions:
\begin{equation}\label{ct8}
(id_Y \otimes {\bf c}_{X,Z}) {\bf a}_{Y,X,Z}({\bf c}_{X,Y}\otimes id_Z)={\bf a}_{Y,Z,X}
{\bf c}_{X,Y\otimes Z} {\bf a}_{X,Y,Z},
\end{equation}
\begin{equation}\label{ct9}
({\bf c}_{X,Z}\otimes id_Y){\bf a}^{-1}_{X,Z,Y}(id_X \otimes {\bf
c}_{Y,Z}) = {\bf a}_{Z,X,Y}^{-1}{\bf c}_{X\otimes Y,Z}{\bf
a}^{-1}_{X,Y,Z}. \end{equation}

If the braiding ${\bf c}$ satisfies ${\bf c}_{X,Y}\cdot
 {\bf c}_{Y,X}=id$ then braided Gr-categories are called {\it symmetric categorical groups}, or {\it Picard categories}.
  Then the relations (\ref{ct9}) and (\ref{ct8}) coincide.

 If $(\mathbb B, \bf c)$ and $(\mathbb B', \bf c')$ are braided Gr-categories, then a
braided Gr-functor   $(F,\Fm): \mathbb B\ri \mathbb B'$ is a
Gr-functor which is compatible with the braidings ${\bf c,c'}$ in the sense that the
following diagram commutes
\[\begin{diagram}\label{bd9}
\node{F(X\tx Y)}\arrow{e,t}{F(\mathbf{c})}
\node{F(Y \tx X)}\\
\node{FX\tx FY}\arrow{e,t}{\mathbf{c}'}\arrow{n,l}{\widetilde{F}}
\node{FY\tx FX}\arrow{n,r}{\widetilde{F}}
\end{diagram}\]

First, let us briefly recall the result on classification of A. Joyal and
R. Street [5].

An {\it abelian 3-cocycle} for $M$ with coefficients in $N$ is a pair
$(h,\eta)$, where $h:M^3\ri N$ is a ``normalized 3-cocycle", satisfying
\begin{equation*}
h(y,z,t)-h(x+y,z,t)+h(x,y+z,t)-h(x,y,z+t)+h(x,y,z)=0,
\end{equation*}
\begin{equation*}
h(x,y,z)-h(y,x,z)+h(y,z,x)+\eta(x,y+z)-\eta(x,y)-\eta(x,z)=0,
\end{equation*}
\begin{equation*}
h(x,y,z)-h(x,z,y)+h(z,x,y)-\eta(x+y,z)+\eta(y,z)+\eta(x,z)=0.
\end{equation*}
For any function $g:M^{2}\rightarrow N$ satisfying
$g(x,0)=g(0,y)=0,$ the {\it coboundary} of $g$ is the abelian
3-cocycle $\partial_{ab}(g)=(h,\eta)$ defined by the equations
\begin{equation*}
h(x,y,z)=g(y,z)-g(x+y,z)+g(x,y+z)-g(x,y),
\end{equation*}
\begin{equation*}
\eta(x,y)=g(x,y)-g(y,x).
\end{equation*}

 A function $\nu:M\rightarrow N$ between abelian groups is called a {\it quadratic map} when the function $M\times M\rightarrow N, (x,y)\mapsto
\nu(x)+\nu(y)-\nu(x+y)$ is bilinear and  $\nu(-x)=\nu(x)$. \\
\indent The {\it trace} of an abelian 3-cocycle $(h,\eta)\in Z^3_{ab}(M,
N)$ is a function
$$t_\eta:M\rightarrow N,\;t_\eta(x)=\eta(x,x).$$
A simple calculation shows that traces are quadratic maps, and
Eilenberg - MacLane [2, 3, 8] proved that the traces determine an
isomorphism
\begin{equation*}
H^{3}_{ab}(M,N)\cong \mathcal{Q}uad(M,N),\;[(h,\eta)]\mapsto
t_{\eta},
\end{equation*}
where $\mathcal{Q}uad(M,N)$ is the abelian group of quadratic maps
from  $M$ to $N$. This result plays a fundamental role in the proof of the
Classification Theorem (Theorem 3.3 [5]).

A. Joyal and R. Street proved that each braided Gr-category $\mathbb
B$ determines a quadratic function $q_\mathbb B:\pi_0\mathbb
B\rightarrow\pi_1\mathbb B$ and let $\mathcal{Q}uad$ be the category
whose objects $(M,N,t)$ are quadratic maps $t:M\ri N$
between abelian groups  $M, N$ and whose morphisms $(\varphi, f): (M, N, t)\ri (M', N', t')$ consist of homomorphisms $\varphi, f$ such that we have a commutative
square
\[\begin{diagram}
\node{M}\arrow{e,t}{\varphi}\arrow{s,l}{t}
\node{M'}\arrow{s,r}{t'}\\
\node{N}\arrow{e,t}{f}
\node{N'}
\end{diagram}\]
Let $\mathcal{BCG}$ denote the category whose objects are braided categorical
groups and whose morphisms are braided monoidal functors. 

\noindent {\bf Theorem 4.1.} (Theorem 3.3 [5]) {\it The functor
\[\begin{matrix}T: &\mathcal{BCG}& \ri& \mathcal{Q}uad,\\
 \;&\mathbb B &\mapsto &(\pi_0\mathbb B, \pi_1\mathbb B, q_{\mathbb
 B})\end{matrix}\]
 has the following properties:}\\
i) {\it For each object $Q$ of $\mathcal{Q}uad$, there exists an object
$\mathbb B$ of
 $\mathcal{BCG}$ with an isomorphism  $T(\mathbb B)\cong Q;$}\\
ii) {\it For any isomorphism $\rho: T(\mathbb B)\stackrel{\sim}{\rightarrow}T(\mathbb B')$, there is an equivalence $F:\mathbb B\ri \mathbb B'$ such that $T(F)=\rho$; and}\\
iii) $T(F)$ {\it is an isomorphism  if and only if $F$ is an
equivalence.}

Now, we state the solution to the classification problem of
braided Gr-categories by the same technique performed for Gr-categories.

If $\mathbb B$ is a braided Gr-category with the braiding $\mathbf c$,
then $\pi_0\mathbb B$ is an abelian group and  acts trivially on
$\pi_1\mathbb B$. Then the reduced Gr-category $S_{\mathbb B}$
becomes a braided Gr-category, with the induced braiding
$\mathbf{c}^\ast=(\bullet, \eta)$ given by the following commutative
diagram:
\[\begin{diagram}
\node{X_r\tx X_s}\arrow{e,t}{i_{X_r\tx X_s}}\arrow{s,l}{\mathbf{c}}
\node{X_{rs}}\arrow{s,r}{\gamma_{_{X_{rs}}}(\eta(r,s))}\\
\node{X_s\tx X_r}\arrow{e,t}{i_{X_s\tx X_r}} \node{X_{sr}}
\end{diagram}\]
Moreover, $(H, \widetilde{H})$ and $(G,\Gm)$ defined by (\ref{ct1}) and (\ref{ct2}) are then braided
Gr-equivalences.

Therefore, each pair $(h,\eta) $ of associativity  and braiding
constraints of $S_{\mathbb B}$ is an abelian 3-cocycle, and $\mathbb
B$ determines uniquely an element $[(h,\eta)]\in
H^3_{ab}(\pi_0\mathbb B,\pi_1\mathbb B)$.

It follows from Theorem 2.5 that

\noindent {\bf Corollary 4.2.} {\it Each braided  Gr-functor $(F,\Fm): \mathbb{S}\ri \mathbb{S}'$ is a
triplet $(\varphi,f, g)$, where}
$$\varphi^\ast (h',\eta')-f_\ast (h,\eta)=\partial_{ab}(g).$$

\noindent Let $\bf H^{3}_{BGr}$ denote the category  whose  objects are
triplets $(M, N,[(h,\eta)])$, where $[(h,\eta)]\in H^3_{ab}(M, N)$.
 A morphism $(\varphi,f):(M, N,[(h,\eta)])\ri (M', N',[(h',\eta')])$
in $\bf H^{3}_{BGr}$  is a pair $(\varphi,f)$ such that there is a
function  $g:M^2\ri N'$ making $(\varphi,f,g):(M, N,(h,\eta))\ri (M', N',(h',\eta'))$ become a braided monoidal functor, i.e.,
$[\varphi^\ast (h',\eta')]=[f_\ast(h,\eta)]\in H^3_{ab}(M, N')$.

We write
$$\text{Hom}_{(\varphi, f)}^{Br}[\mathbb B,\mathbb B']$$
\noindent for the set of homotopy classes of braided Gr-functors
$\mathbb B\ri\mathbb B'$ inducing the same pair $(\varphi, f)$.

 Now,
Corollary 4.2
and the proofs of Theorem 3.1, Theorem 3.2 with some appropriate
modifications lead to the following theorem

\noindent {\bf Theorem 4.3.} (The Classification Theorem) {\it There is a classifying functor

\[\begin{matrix}
d:&\mathcal{BCG}&\ri&\bf H^{3}_{BGr},\\
&\mathbb{B}&\mapsto&(\pi_0\mathbb{B},\pi_1\mathbb{B},[(h,\eta)_\mathbb{B}]),\\
&(F,\widetilde{F})&\mapsto&(F_0,F_1)
\end{matrix}\]
which has the following properties:}\\
i) $dF$ {\it is an isomorphism if and only if $F$ is an
 equivalence,}\\
ii) $d$ {\it is a surjection on  objects,}\\
iii) $d$ {\it is full, but not faithful. For
$(\varphi,f):d\mathbb{B}\rightarrow d\mathbb{B}'$, we have}
$$\text{Hom}^{Br}_{(\varphi,f)}[\mathbb{B},\mathbb{B}']\cong
H^2(\pi_0\mathbb{B},\pi_1\mathbb{B}').$$

We write
$$\mathcal{BCG}[M, N]$$
for  the set of equivalence classes of braided Gr-categories whose
pre-sticks are of the type   $(M, N)$. By Corollary 4.2, we can prove
the following proposition

\noindent {\bf Theorem 4.4.} {\it There exists a bijection}
\begin{eqnarray*}
	\Gamma:\mathcal{BCG}[M, N]&\ri &H^{3}_{ab}(M, N),\\
{[\mathbb B]}&\mapsto & q^{-1}_\ast p^\ast  [(h,\eta)_{\mathbb B}].
\end{eqnarray*}
\\
{\bf 5\ \ Gr-category of an abstract kernel}\\
\\
The notion of \emph{abstract kernel} was introduced in [9]. It is a
triplet $(\Pi, G,\psi)$, where $\psi:\Pi\rightarrow$ Aut$G/$In$G$ is
a group homomorphism. In this section, we describe the
Gr-category structure of an abstract kernel and apply it to make the
constraints of a Gr-category be strict.
The operation of $G$ is denoted by +. 
 The {\it center} of $G,$ denoted by $ZG,$ consists of elements $c\in G$ such that $c+a=a+c$ for all $a\in G.$

Let us recall that the obstruction of $(\Pi, G,\psi)$ is an element $[k]\in
H^{3}(\Pi,ZG)$, defined as follows. For each $x\in \Pi$, choose
$\varphi(x)\in \psi(x)$ such that $\varphi(1)= id_G.$ Then, there is
a function $f:\Pi^2\rightarrow G$ satisfying
 \begin{equation}\label{ct10}
 \varphi(x)\varphi(y)=\mu_{f(x,y)}\varphi(xy).
 \end{equation}
The pair $(\varphi,f)$ therefore induces an element  $k\in Z^{3}(\Pi,ZG)$ by
the relation
\begin{equation}\label{ct11}
\varphi(x)[f(y,z)]+f(x,yz)=k(x,y,z)+f(x,y)+f(xy,z).
\end{equation}
For each group $G$, we can construct a category, denoted by ${\bf Aut}_G$,
whose objects are elements of the group of automorphisms Aut$G$. For
two elements $\alpha,\beta$ of Aut$G,$ we write
$$\text{Hom}(\alpha,\beta)=\{c\in G|\alpha=\mu_c\circ\beta\},$$
where $\mu_c$ is an inner-automorphism  $G,$ induced by $c\in G.$
For two morphisms $c: \alpha\ri \beta;\ d:\beta\ri\gamma$ in ${\bf
Aut}_G$, the composition is defined by  $d\circ c=c+d$ (the addition
in $G$).

The category ${\bf Aut}_G$ is a strict Gr-category with the tensor
product defined by   $\alpha\otimes\beta=\alpha\circ\beta$, and
\begin{equation}\label{ct12}
(\alpha\stackrel{c}{\rightarrow}\alpha')\otimes(\beta\stackrel{d}{\rightarrow}\beta')=
\alpha\otimes\beta\stackrel{c+\alpha'd}{\rightarrow}\alpha'\otimes\beta'.
\end{equation}

The following proposition describes the reduced Gr-category of the
Gr-category of an abstract kernel.

\noindent {\bf Proposition 5.1.} {\it Let $(\Pi,G,\psi)$ be an abstract kernel with $[k]\in H^3(\Pi,ZG)$ be
its obstruction. Let the reduced
Gr-category of the strict one  ${\bf Aut}_G$ be $S_{{\bf Aut}_G}=(\Pi',C,h).$  Then}\\
i) $\Pi'=\pi_0({\bf Aut}_G)=$ Aut$G/$In$G, C=\pi_1({\bf Aut}_G)=ZG,$\\
ii) {\it $\psi^* h$ belongs to the cohomology class of  $k.$}

\begin{proof}
i) It follows from the definition of the category ${\bf Aut}_G$ and
the reduced Gr-category.

\noindent ii) Let $(H, \widetilde{H})$ be a canonical Gr-equivalence from
$\mathbb S$ to ${\bf Aut}_G$. Then, the following diagram
\begin{equation}\label{ct13}
\begin{CD}
Hr(Hs Ht )@>id\tx\widetilde{H}>> HrH(st)@> \widetilde{H}>>H(r(st)) \\
@|                       @.                        @VV H(\bullet,h(r,s,t))V \\
(Hr Hs) Ht @>\widetilde{H}(r,s)\tx id >> H(rs) Ht @> \widetilde{H}>>
H((rs)t)
\end{CD}
\end{equation}
commutes for all $r, s, t\in \Pi'$. Since ${\bf Aut}_G$ is a strict
Gr-category, we have
$$\gamma_\alpha(u)=u, \ \forall \alpha \in \mathbf{Aut}_G,\ \forall u \in ZG=C.$$
Associating  with the definition of $H,$ we obtain $H(\bullet,c)=c,
\ \forall c \in C$. From the commutativity of the diagram  (\ref{ct13}) and
the relation (\ref{ct12}), we have
\begin{equation}\label{ct14}
Hr[g(s,t)]+g(r,st)=g(r,s)+g(rs,t)-h(r,s,t)
\end{equation}
where $g=g_H:\Pi'\times\Pi'\ri G $ is associated with $\widetilde{F}$.\\
\indent For the abstract kernel $(\Pi, G, \psi),$ choose a
function $\varphi=H.\psi:\Pi\ri$ Aut$(G)$. Clearly,
$\varphi(1)=id_G$. Moreover, since
$$\widetilde H_{\psi(x), \psi(y)}:H\psi(x)H\psi(y)\ri H\psi(xy)$$
is a morphism in ${\bf Aut}_G,$  for all $x,y\in\Pi$  we have
$$\varphi(x)\varphi(y)=H\psi(x)H\psi(y)=\mu_{f(x,y)}H\psi(xy)=\mu_{f(x,y)}\varphi(xy),$$
where $f(x,y)=\widetilde H_{\psi(x),\psi(y)}$. Thus, the pair
$(\varphi,f)$ is a factor set of the abstract kernel $(\Pi,G,\psi).$
It induces an {\it obstruction} $k(x,y,z)\in Z^{3}(\Pi,ZG)$
satisfying (\ref{ct11}).
Now, for  $r=\psi(x),\ s=\psi(y),\ t=\psi(z)$, the equation  (\ref{ct14})
becomes
$$\varphi(x)[f(y,z)]+f(x,yz)=+f(x,y)+f(xy,z)-(\psi^* h)(x,y,z).$$
In comparison with (\ref{ct11}), $[\psi^* h]= [k].$
\end{proof}

\noindent We now use Proposition 5.1  and the Theorem on the
realization of the obstruction in the group extension problem to
prove Theorem 5.3. First, we need the following lemma

\noindent {\bf Lemma 5.2.} {\it Let $\mathbb H$ be a strict Gr-category and $S_{\mathbb H}=(\Pi, C,
 h)$ be its reduced Gr-category. Then, for each group
homomorphism  $\psi:\Pi'\ri \Pi$, there exists a strict Gr-category
$\mathbb G$ which is  Gr-equivalent to the Gr-category  $\mathbb
J=(\Pi', C, h')$, where  $C$ is regarded as a $\Pi'$-module with an
operator $xc=\psi(x)c,$ and $h'$ belongs to the same cohomology
class as $\psi^* h.$}
\begin{proof}
We construct a strict Gr-category $\mathbb G$ as follows
\begin{gather*}
\text{Ob}(\mathbb G)=\{(x, X)| \ x\in \Pi', X\in \psi(x)\}, \\
\text{Hom} ((x, X), (x, Y))=\{x\}\times \text{Hom}_{\mathbb H}(X,Y).
\end{gather*}
The tensor products on objects and morphisms of $\mathbb G$ are
defined by
\begin{gather*}
(x, X)\tx (y,Y)=(xy, X\tx Y),\\
(x,u)\tx (y,v)=(xy, u\tx v).
\end{gather*}
The unit object of  $\mathbb G$ is $(1, I)$ where  $I$ is the unit
object of $\mathbb H$. One can verify that   $\mathbb G $ is
a strict Gr-category.  Moreover, we have isomorphisms
\[\begin{aligned}
\lambda:\pi_0\mathbb G&\ri \Pi',\\
 [(x,X)]&\mapsto x,
\end{aligned}\qquad\qquad
\begin{aligned}
\ f:\pi_1\mathbb G&\ri \pi_1\mathbb H=C,\\
(1,c)&\mapsto c,
\end{aligned}\]
\noindent and a  Gr-functor $(F,\Fm):\mathbb G \rightarrow \mathbb
H$ given by
$$F(x,X)=X,\ \ F(x,u)=u, \ \  \Fm=id.$$
 Let  $S_F=(\phi,\widetilde{\phi} ):S_{\mathbb G}\rightarrow S_{\mathbb H}$ be a
Gr-functor induced by $(F,\Fm)$ between the reduced categories.
Then, we have
\begin{gather*}
\phi(x,X)=F_0(x,X)= [F(x,X)]= [X]=\psi(x),\\
\phi(1,u)=F_1(1,u)=\gamma_{F(1,I)}F(1,u)=\gamma_I(u)=u,
\end{gather*}
where $u$ is a morphism in $\mathbb G$. This means $F_0=\psi
\lambda$ and $F_1=f$, or $S_F$ is
a functor of the type  $(\psi\lambda, f).$\\
\indent Now if $ h_{\mathbb G}$ is the associativity constraint of
$S_{\mathbb G}$.
By Theorem 2.6, the obstruction of the pair $(\psi\lambda, f)$ must
vanish in $H^{3}(\pi_0\mathbb G, \pi_1\mathbb H)=H^{3}(\pi_0\mathbb
G, C)$, i.e.,
$$ (\psi\lambda)^{\ast} h=f_{\ast} h_{\mathbb G} + \delta \widetilde{\phi}. $$
Now, if we denote $h'=f_{\ast} h_{\mathbb G}$, the pair $J=(\lambda,
f), \widetilde{J}=id$ is a Gr-functor from  $S_\mathbb G$ to
$\mathbb J=(\Pi',C,h').$ Then, the composition
$$\mathbb G \stackrel{(G,\widetilde{G})}{\longrightarrow}S_{\mathbb G}\stackrel{(J,\widetilde{J})}{\longrightarrow}\mathbb J$$
 is a Gr-equivalence from $\mathbb G$ to $\mathbb J=(\Pi',C, h').$

Finally, we prove that $ h'$ belongs to the same cohomology
class as  $\psi^{\ast} h.$ Let $K=(\lambda^{-1}, f^{-1}):(\Pi', C,
h')\rightarrow S_{\mathbb G}$. Then $K$ together with
$\widetilde{K}=id$ is a Gr-functor, and the composition
$$(\phi, \widetilde{\phi})\circ (K, \widetilde{K}):(\Pi', C,
h')\rightarrow S_{\mathbb H}$$ is a  Gr-functor making the following
diagram commute
\[
\begin{diagram}
\node[1]{S_{\mathbb G}}
 \arrow[2]{e,t}{\phi}
\node[2]{S_{\mathbb H}}
\\
\node[2]{\mathbb J=(\Pi', C, h')}\arrow[1]{nw,r,2}{K}
\arrow[1]{ne,r,2}{\phi\circ K}
\end{diagram}
\]
 Clearly, $\phi\circ K$ is a Gr-functor of the type $(\psi, id)$ and therefore its obstruction vanishes. By (\ref{ct4}), we have $\psi^*h-h'=\partial g,$ i.e., $[h']=[\psi^*h].$

\end{proof}
\noindent {\bf Theorem 5.3.} {\it Each Gr-category is  Gr-equivalent to a strict one.}
\begin{proof}
 Let $\mathbb C$ be a Gr-category
whose reduced Gr-category is $S_{\mathbb C}=(\Pi', C, k)$. By the
theorem on the realization of the obstruction  (Theorem 9.2 Section
IV [9]), the realization of 3-cocycle $k\in H^{3}(\Pi',C)$  is the
group $G$ with the center  $ZG=C$ and group homomorphism
$\psi:\Pi'\ri$ Aut$G/$In$G$ such that $\psi$ induces a
$\Pi'$-module structure on $C$ and
 the obstruction of the  abstract kernel $(\Pi', G, \psi)$ is $k$. By Proposition 5.1,
  $S_{{\bf Aut}_G}=$(Aut$G/$In$G, C,  h)$ is the reduced Gr-category of the strict Gr-category
    ${\bf Aut}_G$, where  $[\psi^* h]= [k].$

Using Lemma 5.2 for $\mathbb H={\bf Aut}_G$, the homomorphism
$\psi:\Pi'\ri$ Aut$G/$In$G$ defines a strict Gr-category $\mathbb
G$, which is Gr-equivalent to the strict Gr-category $\mathbb J=(\Pi', C,
h').$ The $\Pi'$-module structures of  $C$ on $S_{\mathbb C}$ and on
$\mathbb J$ coincide. Moreover,  $[\psi^* h]= [h'].$ It follows that
$[h']= [k].$ So there is a function $g:\Pi'\times \Pi'\ri C$ such
that $h'-k=\partial g$. Then, by Theorem 2.6, $$ (K, \widetilde{K})=
(id_{\Pi'}, id_C, g):S_\mathbb C\ri \mathbb J$$ \noindent is a
Gr-equivalence. Therefore, $\mathbb C$  is equivalent to the strict
Gr-category  $\mathbb G$.
\end{proof}
 Readers can see a different proof of Theorem 5.3 in [15].\\
 \\
{\bf 6\ \ Gr-functors and the group extension problem}\\
\\
In this section, we apply Theorem 3.1 to obtain Schreier classical
Theorem on group extensions.

\noindent {\bf Theorem 6.1.} {\it Let $G$ and $\Pi$ be groups. Then}\\
i) {\it There is a canonical partition}
$$\text{Ext}(\Pi,G)=\underset{\psi}\coprod
\text{Ext}(\Pi,G,\psi),$$
\noindent {\it where, for each morphism} $\psi:\Pi\ri{\bf Aut}G/\text{In}G,$ Ext$(\Pi,G,\psi)$
{\it is the set of equivalence classes of group extensions $E:G\ri
B\ri \Pi$ of
$G$ by $\Pi$ which induce $\psi$.}\\
ii) {\it Each abstract kernel $(\Pi,G,\psi)$ determines a (normalized)
third-dimensional cohomology class } Obs$(\Pi,G,\psi)\in H^3(\Pi,ZG)$
{\it (with respect to the $\Pi$-module structure on $ZG$ obtained
via
 $\psi$), called the obstruction of  $(\Pi,G,\psi)$. The abstract
 kernel has extensions if and  only if its obstruction vanishes.
 Then, there is a bijection}
$$\text{Ext}(\Pi,G)\leftrightarrow H^2_{\psi}(\Pi,ZG).$$
As below, each factor set $(\varphi,f)$ of a group extension can
be lifted to a Gr-functor $F:\text{Dis}\Pi\rightarrow{\bf Aut}G$,
when Dis$\Pi$ is regarded as a Gr-category of the type $(\Pi, 0, 0)$, and
therefore we can classify group extensions by means of Gr-functors.

We write $\text{Hom}_\psi[$Dis$\Pi,{\bf Aut}G]$ for the set of
homotopy class of Gr-functors from Dis$\Pi$ to {\bf Aut}G inducing the pair of homomorphisms $(\psi,0)$ and
$\text{Ext}_\psi(\Pi,G)$ for the set of equivalence classes of group
extensions of $G$ by $\Pi$ inducing $\psi$, we have

\noindent {\bf Theorem 6.2.} {\it There exists a bijection}
$$\Delta:\text{Hom}_\psi[\text{Dis}\Pi,{\bf Aut}_G]\ri\text{Ext}_\psi(\Pi,G).$$
\begin{proof}

{\it Step 1: The construction of the group extension $E_F$ of  $G$ by
$\Pi,$ induced by Gr-functor $F.$}

Let $(F,\Fm):\text{Dis}\Pi\rightarrow{\bf Aut}G$ be a Gr-functor.
Then, $\Fm_{x,y}=f(x,y)$ is a function from $\Pi^2$ to $G$ such
that
\begin{equation}\label{ct15}
Fx\circ Fy=\mu_{f(x,y)}\circ F{xy}.
\end{equation}
The compatibility of  $(F,\Fm)$ with the unit and associativity constraints, respectively, implies
\begin{equation}\label{ct16}
Fx[f(y,z)]+f(x,yz)=f(x,y)+f(xy,z),
\end{equation}
\begin{equation}\label{ct17}
f(x,1)=f(1,y)=0.
\end{equation}
Set $B_F=\{(a,x)|a\in G,x\in\Pi\}$ and the operation
$$(a,x)+(b,y)=(a+Fx(b)+f(x,y),xy).$$
Then, $B_F$ is an extension of $G$ by $\Pi,$
$$E_F:\;\;0\ri G\stackrel{i}{\ri}B_F\stackrel{p}{\ri}\Pi\ri 1,$$
where $i(a)=(a,1),p(a,x)=x.$ The relations (\ref{ct15}), (\ref{ct16}) imply the
associativity  of the operation in $B_F.$ Indeed,
 the unit of the addition  in $B_F$ is $(0,1),$
 the opposite element $(a,x)\in B_F$ is $(b,x^{-1})\in B_F,$ where $b$ is an element such that $Fx(b)=-a+f(x,x^{-1}).$\\
 The {\it conjugation} homomorphism $\psi:\Pi\rightarrow$ Aut$G/$In$G$ is determined by $\psi(x)=[\mu_{(0,x)}]$. By a simple
 calculation, we have  $\mu_{(0,x)}(a,1)=(Fx(a),1)$. Let $G$ and its image $iG$ be identical, we obtain $\psi(x)=[Fx].$

\noindent {\it Step 2: $F$ and $F'$ are homotopic if and only if $E_F$ and
$E_{F'}$ are congruent.}

Let $F,F':\text{Dis}\Pi\ri{\bf Aut}_G$ be Gr-functors and
$\alpha:F\ri F'$ be a homotopy. Then, by the definition of
Gr-morphisms, the following diagram commutes
\[
\begin{diagram}
\node{Fx\otimes Fy}\arrow{e,t}{\Fm}
\arrow{s,l}{\alpha_x\otimes\alpha_y}
\node{F(xy)}
\arrow{s,r}{\alpha_{xy}}\\
\node{F'x\otimes F'y}\arrow{e,t}{\Fm'}
\node{F'(xy)}
\end{diagram}
\]
That is,
$$
\Fm_{x,y}+\alpha_{xy}=\alpha_x\otimes\alpha_y+\Fm'_{x,y},
$$
or
\begin{equation}\label{ct18}
f(x,y)+\alpha_{xy}=\alpha_x+F'x(\alpha_y)+f'(x,y).
\end{equation}
Now, we write
\begin{eqnarray}
    \beta:B_F&\ri&B_{F'},\nonumber\\
    (a,x)&\mapsto&(a+\alpha_x,x).\nonumber
\end{eqnarray}
Note that $Fx=\mu_{\alpha_x}\circ F'x$, and by (\ref{ct15}) one can see that
$\beta$ is a homomorphism. Moreover, it is an isomorphism making the
following diagram commute, i.e.,  $E_F$ and $E_{F'}$ are congruent.
\[\divide \dgARROWLENGTH by 2
\begin{diagram}
\node{E_F:\;\;\;0} \arrow{e,t}{} \node{G} \arrow{e,t}{i}
\arrow{s,r}{id} \node{B_F} \arrow{e,t}{p}\arrow{s,r}{\beta}
\node{\Pi} \arrow{e,t}{}\arrow{s,r}{id}\node{1}\\
\node{E_{F'}:\;\;\;0} \arrow{e,t}{} \node{G} \arrow{e,t}{i'}
\node{B_{F'}} \arrow{e,t}{p'} \node{\Pi} \arrow{e,t}{} \node{1}
\end{diagram}
\]
The conversion of the proposition can be obtained as we see by retracing our steps.

\noindent {\it Step 3:  $\Delta$ is a surjection}.

Suppose that the group extension
$$E:0\ri G\stackrel{i}{\ri}B\stackrel{p}{\ri}\Pi\ri 1,$$
associates with the homomorphism  $\psi:\Pi\rightarrow$ Aut$G/$In$G$.
Select a ``representative'' $u_x,x\in \Pi,$ in $B$, that is
$p(u_x)=x$. In particular, choose $u_1=0$. Then, the elements of $B$ can
be written  uniquely as $a+u_x$, for $a\in G, x\in \Pi$, and
$$u_x+a=\mu_{u_x}(a)+u_x.$$

The sum $u_x+u_y$ must be  in the same coset as $u_{xy}$, so there
are unique elements $f(x,y)\in G$ such that
\begin{equation*}
u_x+u_y=f(x,y)+u_{xy}.
\end{equation*}

The function $f$ is called a {\it factor set} of the extension $E$.
Thus, it satisfies the relations
\begin{equation}\label{ct19}
\mu_{u_x}[f(y,z)]+f(x,yz)=f(x,y)+f(xy,z),\;\;\;x,y,z\in\Pi.
\end{equation}
\begin{equation}\label{ct20}
f(x,1)=f(1,y)=0.
\end{equation}

We construct a Gr-functor $F=(F,\Fm)$: Dis$\Pi\ri{\bf Aut}_G$ as
follows: $Fx=\mu_{u_x}$, $\Fm_{x,y}=f(x,y).$

Clearly, the relations (\ref{ct19}), (\ref{ct20}) show that $(F,\Fm)$ is
a monoidal functor between Gr-categories.
\end{proof}
We now prove Theorem 6.1.

Let $(\Pi,G,\psi)$ be an abstract kernel. For each $x\in\Pi$, choose $\varphi(x)\in\psi(x)$  such that $\varphi(1)=id_G$. The
family of $\varphi(x)$ induces a function $f:\Pi^2\ri G$  satisfying
the relation (\ref{ct10}). The pair $(\varphi,f)$ induces an obstruction
$k\in Z^3(\Pi,ZG)$ by the relation (\ref{ct11}). Write $F(x)=\varphi(x)$,
we obtain a functor Dis$\Pi\ri {\bf Aut}_G$.

Let $\mathbb S=$ (Aut$G/$In$G,ZG, h)$ be the reduced Gr-category of
 ${\bf Aut}_G$. Then $F$ induces the pair of group homomorphisms
 $(\psi,0):(\Pi,0)\ri ($Aut$G/$In$G,ZG)$ and by the relation (\ref{ct4}) an obstruction of the
 functor $F$ is $\psi^*h$. By Proposition 5.1, $[\psi^*h]=[k]$, i.e., the obstruction   of the
abstract kernel $(\Pi,G,\psi)$ and the obstruction of the functor
$F$ coincide. Then, by Theorem 2.6, $(\Pi,G,\psi)$ has extensions if
and only if its obstruction vanishes.

According to Theorem 3.1, there is a bijection
$$\text{Hom}_{(\psi,0)}\text{[Dis}\Pi, {\bf Aut}_G]\leftrightarrow H^2(\Pi,ZG),$$
since $\pi_0($Dis$\Pi)=\Pi, \pi_1({\bf Aut}_G)=ZG$. Combination with
Theorem 6.1 yields:

$$\text{Ext}(\Pi,G)\leftrightarrow H^2(\Pi,ZG).$$
This completes the proof.

\begin{center}

\end{center}
\end{document}